\documentclass{elsart}

\usepackage{amssymb}

\newtheorem{theorem}{Theorem}[section]
\newtheorem{lemma}[theorem]{Lemma}

\newtheorem{corollary}[theorem]{Corollary}

\newtheorem{example}[theorem]{Example}

\newcommand \Q {\mathbb Q}

\newcommand \N {\mathbb N}
\newcommand \bpig {{\textup B}_{\pi}(G)}

\newcommand \ipi {{\textup I}_{\pi}}
\newcommand \ipig {{\textup I}_{\pi}(G)}
\newcommand \npig {{\textup N}_{\pi}(G)}
\newcommand \npi {{\textup N}_{\pi}}
\newcommand \irr {\textup{Irr}}
\newcommand \irrg {\textup{Irr}(G)}
\newcommand \ibr {{\textup{IBr}}_p}
\newcommand \ibrg {{\textup{IBr}}_p(G)}
\newcommand \sbs {\subseteq}
\newcommand \nnq {{\bf{N}}_N(Q)}
\newcommand \nmq {{\bf{N}}_M(Q)}
\newcommand \nnmq {{\bf{N}}_{NM}(Q)}
\newcommand \nwq {{\bf{N}}_W(Q)}
\newcommand \nuq {{\bf{N}}_U(Q)}
\newcommand \ngq {{\bf{N}}_G(Q)}
\newcommand \wt {\widetilde}
\newcommand \dvd {\hbox {\Big|}}
\newcommand \ndvd {\hbox {\big/}\kern-5pt\dvd}
\newcommand \nrml {\lhd}
\def \< {\langle}
\def \> {\rangle}

\begin{document}

\begin{frontmatter}



\title{Bounds on the number of lifts of a Brauer character in a p-solvable group}


\author{James P. Cossey}

\address{Department of Mathematics, University of Arizona, Tucson, AZ 85721}
\ead{cossey@math.arizona.edu} 

\begin{abstract}
  The Fong-Swan theorem \cite{blocks} shows that for a $p$-solvable group
$G$ and Brauer character $\varphi \in \ibrg$, there is an ordinary character $\chi \in \irrg$ such that
$\chi^0 = \varphi$, where $^0$ denotes restriction to the $p$-regular
elements of $G$.  This still holds in the generality of $\pi$-separable
groups \cite{bpi}, where $\ibrg$ is replaced by $\ipig$.  For $\varphi  
\in \ipig$, let $L_{\varphi} = \{ \chi \in \irrg \mid  \chi^0 = \varphi \}$.
In this paper we give a lower bound for the size of $L_{\varphi}$ in
terms of the structure of the normal nucleus of $\varphi$ and, if $G$ is
assumed to be odd and $\pi = \{p' \}$, we give an
upper bound for $L_{\varphi}$ in terms of the vertex subgroup for $\varphi$.
\end{abstract}

\end{frontmatter}

\section{Introduction}

The Fong-Swan theorem (see \cite{blocks})
asserts that if $G$ is a finite $p$-solvable group, and $\varphi$ is an
irreducible  Brauer character of $G$, then there must exist an ordinary character $\chi \in \irrg$ such that $\chi^0 = \varphi$, where $^0$ denotes restriction
to the set of $p$-regular elements of $G$.  Such a character is called a lift of $\varphi$.  Moreover, the set $\ibrg$ of the irreducible Brauer characters of $G$ forms a basis of the
space of class functions on the $p$-regular elements of $G$ and $G$ has the
property that if $\chi \in \irrg$, then $\chi^0$ is an $\N$-linear
combination of elements of $\ibrg$.

Let $\pi$ be any set of primes.  Isaacs in \cite{bpi} has generalized the above to the case where $G$ is
$\pi$-separable, and the set $\ibrg$ is replaced by the set
$\ipig$, the irreducible $\pi$-partial characters of $G$.  Again, if $\chi                
\in \irrg$, it is shown that $\chi^0$ is an $\N$-linear combination of
elements of $\ipig$, where here $^0$ denotes
restriction to the set of $\pi$-elements.  Isaacs then defines
a set $\bpig \sbs \irrg$ of lifts of $\ipig$ such that for
each $\varphi \in \ipig$, there is a unique character $\chi \in \bpig$
such that $\chi^0 = \varphi$.  There are other ways, however, to
construct sets of lifts of the $\ipig$ characters, as in
\cite{npiconst}, where a canonical set of lifts which we will call $\npig$ is constructed.  Little is known about the set of all lifts of a
particular character $\varphi \in \ipig$.  For instance, using character triples, Laradji in \cite{twobrauer} has given lower bounds for the number of such lifts in the case that $p = 2$.

Let $G$ be a $\pi$-separable group.  For a fixed character $\varphi \in \ipig$, we define $L_{\varphi}$ to be $\{ \chi \in \irrg \mid  \chi^0 = \varphi \}$.  This set will be the focus of study in this paper.  The first main result of this paper,
which follows immediately from the work of Navarro
\cite{newpispecial}, is a lower bound on the size of $L_{\varphi}$ in terms of a certain subgroup $W_{\varphi}$ determined by $\varphi$.

\begin{theorem}  Let $G$ be a $\pi$-separable group, and let $\varphi \in \ipig$.  Then there is a 
subgroup $W_{\varphi}$ and an irreducible $\pi$-partial character $\alpha$ of $W_{\varphi}$ such that $| L_{\varphi}|  \geq | W_{\varphi}:  
W_{\varphi}'| _{\pi'}$, $\alpha^G = \varphi$, and $\alpha(1)$ is a $\pi$-number.
\end{theorem}

We will give an example where this
inequality is strict.  It is not known if there is a better lower
bound.  The main result of this paper is the following, which gives an
upper bound in the case when $G$ has odd order and $\pi = p'$.  We define the  vertex $Q$ of $\varphi$ to be simply a Sylow $p$-subgroup of
$W_{\varphi}$.

\begin{theorem}  Suppose $G$ has odd order, and suppose $\varphi \in                                   
\ibrg$ for some prime p.  If $Q$ is a vertex subgroup for $\varphi$,
then $| L_{\varphi}|  \leq | Q: Q'| $.
\end{theorem}

Combining the above results, we see that if $G$ has odd order and
$\varphi \in \textup{I}_{p'}(G) = \ibrg$, then $| W_{\varphi}: W_{\varphi}'| _p \leq | L_{\varphi}|  \leq |Q:Q'| $.  Thus if the two above indices are equal (for instance if $W_{\varphi}$ is nilpotent), then we get an exact count on the size of $L_{\varphi}$.

\section{Properties of $\pi$-special characters}
In this section we summarize the properties of $\pi$-special and $\pi$-factorable characters of a finite $\pi$-separable group $G$.  All of the results in this section can be found in \cite{gajen} and \cite{manzwolf}.  For the remainder of this paper, assume that $G$ is $\pi$-separable, unless otherwise indicated.

Recall that the order $o(\chi)$ of an irreducible character $\chi$ is the order of the linear character $\lambda = {\textup{det}}(\chi)$, i.e. the order $o(\chi)$ is the smallest positive integer $n$ such that $({\textup{det}}(\chi(x)))^n = 1$ for every element $x$ of $G$.  If $G$ is a $\pi$-separable group, a character $\alpha \in \irrg$ is called
$\pi$-special if (a) $\alpha(1)$ is a $\pi$-number and (b) for every subnormal subgroup 
$S$ of $G$, and every irreducible constituent $\gamma$ of $\alpha_S$, we have that the order of $\gamma$ is a
$\pi$-number.  Recall that $\pi'$ denotes the complement of the set of primes $\pi$.  If
$\alpha$ and $\beta$ are $\pi$-special and $\pi'$-special,
respectively, then one can show that in fact $\chi = \alpha \beta$ is irreducible.  If an irreducible character $\chi$ can be written in the form $\chi = \alpha \beta$ where $\alpha$ is $\pi$-special and $\beta$ is $\pi'$-special, we say that $\chi$ is $\pi$-factorable.  If $\chi = \alpha_1 \beta_1$ is
another such factorization, then $\alpha_1 = \alpha$ and $\beta_1 = \beta$, i.e. the factorization is unique.  We will frequently write that if $\gamma \in \irrg$ is $\pi$-factorable, then $\gamma = \gamma_{\pi} \gamma_{\pi'}$, where $\gamma_{\pi}$ and $\gamma_{\pi'}$ are the $\pi$-special and $\pi'$ special factors of $\gamma$, respectively.  Note that if $\gamma \in \irrg$ is $\pi$-factorable, and if $S \nrml \nrml G$, then all of the constituents of $\gamma_S$ must be $\pi$-factorable.

The $\pi$-special characters of $G$ (denoted by $\mathcal{X}_{\pi}(G)$)
behave well with respect to normal subgroups.  It is clear from the
definition that if $\alpha \in \mathcal{X}_{\pi}(G)$ and $S \nrml \nrml G$ then
all the irreducible constituents of $\alpha_S$ are in
$\mathcal{X}_{\pi}(S)$.  If $N \nrml G$ and $G/N$ is a $\pi'$-group and if $\alpha \in \mathcal{X}_{\pi}(G)$, then $\alpha$ restricts irreducibly to $\alpha_N \in               
\mathcal{X}_{\pi}(N)$.  Moreover, if $\gamma \in \mathcal{X}_{\pi}(N)$ and
$\gamma$ is invariant in $G$, then there is a unique $\pi$-special
extension of $\gamma$ to $G$.  If $M \nrml G$ and $G/M$ is a
$\pi$-group and if $\delta \in \mathcal{X}_{\pi}(M)$, then every irreducible
constituent of $\gamma^G$ is in $\mathcal{X}_{\pi}(G)$.

Finally, we note that $\mathcal{X}_{\pi}(G)$ is related to $\ipig$.
If $\varphi \in \ipig$ and $\varphi(1)$ is a $\pi$-number, then there
necessarily exists a unique character $\alpha \in \mathcal{X}_{\pi}(G)$ such that
$\alpha^0 = \varphi$.  Also, if $\gamma \in \mathcal{X}_{\pi}(G)$, then
$\gamma^0 \in \ipig$.  In fact, if $\gamma \in \mathcal{X}_{\pi}(G)$ and if $H$ is a Hall $\pi$-subgroup of $G$, then $\gamma$ restricts irreducibly to $H$, and $\gamma$ is the unique $\pi$-special extension of $\gamma_H$.

\section{The normal nucleus and vertices}

In this section we summarize Navarro's construction of the normal nucleus of $\chi \in \irrg$, and we define a set $\npig$ of canonical lifts of
$\ipig$.  We then define the vertex pair $(Q, \delta)$ for $\chi               
\in \irrg$.  Proofs of these results can all be found in \cite{npiconst}.

In \cite{npiconst}, Navarro defines, for a $p$-solvable group $G$ and a character $\chi \in \irrg$, the normal nucleus $(W, \gamma)$ of $\chi$.  In
\cite{npiconst}, this normal nucleus is defined only for $p$-solvable
groups, but the same definitions and results go through if $G$ is
assumed only to be $\pi$-separable for some set of primes $\pi$.  We
here briefly review this construction and highlight some of the
properties we need.

The $\pi$-analog of Corollary 2.3 of \cite{npiconst} states that if
$G$ is $\pi$-separable and if $\chi \in \irrg$, then there is a unique normal subgroup $N$ of $G$ maximal with the property that if $\theta \in \irr(N)$ lies under $\chi$, then $\theta$ is $\pi$-factorable.  (Compare this to
Theorem 3.2 of \cite{bpi}.)  We call the pair $(N, \theta)$ a maximal factorable normal pair of $G$ lying under $\chi$.  Let $G_{\theta}$ be the stabilizer of $\theta$ and let $\psi \in \irr(G_{\theta})$ be the Clifford
correspondent for $\chi$ lying above $\theta$.  If $\theta$ is invariant
in $G$, it is shown that $G = N$, so $\theta = \chi$.  In this case,
the normal nucleus of $\chi$ is defined to be the pair $(G, \chi)$,
and note that $\chi = \theta$ is $\pi$-factorable.  If $G_{\theta} < G$,
then the normal nucleus $(W, \gamma)$ of $\chi$ is defined recursively
to be the normal nucleus of $\psi \in \irr(G_{\theta} \mid  \theta)$.
Note, then, that $\gamma$ is $\pi$-factorable, and $\gamma^G = \chi$.
Also, the normal nucleus $(W, \gamma)$ of $\chi$ is uniquely defined
up to conjugacy.

If $G$ is $\pi$-separable and $\chi \in \irrg$ has the normal nucleus $(W, \gamma)$, we define the
vertex pair of $\chi$ to be the pair $(Q, \delta)$, where $Q$ is a Hall
$\pi'$-subgroup of $W$ and $\delta = (\gamma_{\pi'})_Q$, which is
necessarily in $\irr(Q)$ since $\gamma_{\pi'}$ is $\pi'$-special.
Since $(W, \gamma)$ is unique up to conjugacy, then $(Q, \delta)$ is
unique up to conjugacy.  We define ${\textup{vtx}}(\chi) = (Q, \delta)$.  In the case that $G$ is $p$-solvable, Navarro defines the set $\irr(G | Q, \delta)$ as those irreducible characters of $G$ with vertex pair $(Q,               
\delta)$.  It is then shown that as $Q$ runs over all of the $p$-subgroups of $G$, then $\irr(G |  Q, 1_Q)$ is a set of
lifts of $\ibrg$, and each $\varphi \in \ibrg$ has precisely one lift in that set.  In the general case where $G$ is
$\pi$-separable, we introduce a slightly different notation and define
$\npig$ to be the
set of characters of $G$ with a $\pi$-special normal nucleus character, i.e. those irreducible characters of $G$ with a trivial vertex character.  Again,
the set $\npig$ forms a canonical set of lifts of $\ipig$.
Also, the characters of $\npig$ with $\pi$-degree are precisely
the $\pi$-special characters of $G$.

Classically, a vertex subgroup of a character $\varphi \in \ibrg$ of a $p$-solvable group $G$ has been defined as a Sylow $p$-subgroup of any subgroup $U$ of $G$ with the property that there is an irreducible Brauer character $\alpha$ of $U$ with $p'$-degree such that $\alpha^G = \varphi$ (see \cite{blocks}).  This same definition can be extended to a $\pi$-separable group $G$ if $\pi$ plays the role of $p'$.  Note that if $\chi \in \npig$ is a lift of $\varphi \in \ipig$, then if $(W, \gamma)$ is a normal nucleus character of $\chi$, then $\gamma$ is necessarily $\pi$-special, and $\gamma^0 = \alpha \in \ipig$ is such that $\alpha^G = \varphi$ and $\alpha$ has $\pi$-degree.  Thus if $Q$ is a Hall $\pi'$-subgroup of $W$, then $Q$ is a vertex in the classical sense for $\varphi$.  The following theorem, from \cite{isaacsnavarro}, shows essentially that vertex subgroups of $\varphi$ are the same no matter how one defines them.

\begin{theorem}  Suppose $G$ is $\pi$-separable, and let $\varphi \in \ipig$.  Then there exists a subgroup $U \subseteq G$ and a character $\alpha \in \ipi(U)$ such that $\alpha^G = \varphi$ and $\alpha(1)$ is a
$\pi$-number.  If we let $Q$ be a Hall $\pi$-complement in $U$, then
up to $G$-conjugacy, $Q$ is uniquely determined by $\varphi$, and it is
independent of the choice of $U$.

\end{theorem}

The existence part of the above theorem can be established by the construction of $\npig$.  The uniqueness, however, requires some work to prove.

The construction of the set $\npig$ is obviously very similar to the construction of the set $\bpig$ in \cite{bpi}.  One might ask if these two sets are in fact equal.  It is shown in \cite{counterexample} that $\bpig = \npig$ if $|G|$ is odd or if $2 \in \pi$, but they need not be equal if $2 \in \pi'$.

\section{Theorem 1.1 and related results}

In this section we prove some easy preliminary results and then prove Theorem 1.1 of the introduction.

The following argument to prove Theorem 4.1 is essentially the same as Navarro's proof of Theorem C of \cite{newpispecial}, only with the subnormal nucleus being replaced by the normal nucleus.

\begin{theorem}  Suppose that $G$ is a $\pi$-separable group and suppose $\chi \in \npig$ has normal nucleus $(W, \alpha)$.  If $\gamma \in \mathcal{X}_{\pi'}(W)$, then $(\alpha \gamma)^G \in \irrg$.  If $\beta \in \mathcal{X}_{\pi'}(W)$, then $(\alpha \gamma)^G = (\alpha \beta)^G$ if and only if $\gamma = \beta$.
\end{theorem}

We prove the first statement as a separate lemma.

\begin{lemma}:  Let $\chi \in \npig$ have normal nucleus $(W, \alpha)$.  If $\gamma \in \mathcal{X}_{\pi'}(W)$, then $(\alpha\gamma)^G \in \irrg$.
\end{lemma}

\begin{pf}  We prove the lemma by induction on $|G|$.  Let $(N, \theta)$ be a
maximal factorable normal pair lying under $\chi$.  If $\theta$ is invariant
in $G$, then $N = W = G$ and $\chi$ is $\pi$-special, and thus in this case $\chi \gamma \in \irrg$.

Otherwise we can assume that $N < G$ and thus $G_{\theta} < G$.  For a character $\gamma \in \mathcal{X}_{\pi'}(W)$, let $\beta \in \irr(N)$ be a constituent of $\gamma_N$, so
$\beta \in \mathcal{X}_{\pi'}(N)$.  We claim that $(N, \theta \beta)$
is a maximal factorable normal pair in $G$.  Since $\theta$ is $\pi$-special and $(N, \theta)$ is a
maximal factorable pair in $G$, then $\mathbb{O}_{\pi}(G/N) = 1$.  If $(M,     \rho)$ was a factorable normal pair over $(N, \theta \beta)$ with
$M/N$ a $\pi'$-group, then $\theta$ must extend to $\rho_{\pi}$,
contradicting the maximality of $(N, \theta)$.  Thus $(N, \theta \beta)$ is a maximal factorable pair in $G$.

Suppose $\eta \in \npi(G_{\theta})$ is the Clifford correspondent for $\chi$ lying over $\theta$.  By the inductive hypothesis, we have that $(\alpha \gamma)^{G_{\theta}} \in \irr(G_{\theta})$, since $(W, \alpha)$ is also a normal nucleus for $\eta$.  Now $(\alpha \gamma)^{G_{\theta}}$ lies over $\theta \beta$, and $G_{\theta} \cap G_{\beta} = G_{\theta \beta}$ by the
uniqueness of factorization, and so there is a character $\xi \in \irr(G_{\theta \beta} \mid  \theta \beta)$ that induces to $(\alpha \gamma)^{G_{\theta}}$.  But $\xi$ induces irreducibly to $G$, so $\alpha \gamma$ must induce irreducibly to $G$.
\end{pf}

\begin{lemma}  Suppose the characters $\chi$ and $\psi$ in $\npig$ have normal nuclei $(W, \alpha)$ and $(U, \gamma)$, respectively.  If $\beta \in \mathcal{X}_{\pi'}(W)$ and $\delta \in \mathcal{X}_{\pi'}(U)$ are such
that $(\alpha \beta)^G = (\gamma \delta)^G$, then there is an element $x \in G$
such that $(W, \alpha \beta)^x = (U, \gamma \delta)$.
\end{lemma}

\begin{pf}  Let $\eta = (\alpha \beta)^G = (\gamma \delta)^G$.
If $(N_1, \theta_1)$ and $(N_2, \theta_2)$ are maximal factorable
normal pairs lying under $\chi$ and $\psi$, respectively, then note that $\theta_i$ is
$\pi$-special.  By the same argument as in the above lemma, if
$\epsilon_1 \in \mathcal{X}_{\pi'}(N_1)$ lies under $\beta$ and
$\epsilon_2 \in \mathcal{X}_{\pi'}(N_2)$ lies under $\delta$, then
$(N_1, \theta_1 \epsilon_1)$ and $(N_2, \theta_2 \epsilon_2)$ are maximal factorable normal pairs lying 
under $\eta$.  Thus $N_1 = N_2$ and $\theta_1 \epsilon_1$ is conjugate
to $\theta_2 \epsilon_2$ and thus by the uniqueness of
$\pi$-factorability, $\theta_2$ is conjugate to $\theta_1$ and
$\epsilon_2$ is conjugate to $\epsilon_1$.  Thus if $\chi$ is
$\pi$-special then $\psi$ is $\pi$-special, and in this case
the statement of the lemma follows immediately.

Otherwise we may assume that neither $\chi$ nor $\psi$ is $\pi$-special, and replacing by $G$-conjugates if necessary, we may also assume
that $\theta_1 = \theta_2$ and $\epsilon_1 = \epsilon_2$ and
$G_{\theta} < G$, where $\theta = \theta_1$ and $\epsilon =  \epsilon_1$.  Let $\xi \in \irr(G_{\theta \epsilon} \mid  \theta \epsilon\
)$ be the Clifford correspondent for $\eta$.  Let $\mu_1 \in                   \npi(G_{\theta}|  \theta)$ be the Clifford correspondent for $\chi$ and
$\mu_2 \in \npi(G_{\theta}|  \theta)$ be the Clifford correspondent for
$\psi$.  Since $\alpha \beta$ and $\gamma \delta$ both lie over
$\theta \epsilon$ and induce to $\eta \in \irrg$, then they must
induce to the same character $\xi^{G_{\theta}} = \nu \in \irr(G_{\theta})$ since $G_{\theta \epsilon} \subseteq G_{\theta}$.  Thus the induction hypothesis
applied to $G_{\theta}$ and the character $\nu =  (\alpha                      
\beta)^{G_{\theta}} = (\gamma \delta)^{G_{\theta}}$ and the
$\npi(G_{\theta})$ characters $\mu_1$ and $\mu_2$ gives that $(W, \            \alpha \beta)$ is $G_{\theta}$-conjugate to $(U, \gamma \delta)$.
\end{pf}

We are now ready to prove Theorem 4.1.

\begin{pf} (of Theorem 4.1)  We have proven the first statement,
and one direction of the second statement is obvious.  To complete the
proof, apply the above lemma to the case $\chi = \psi$ and the normal nucleus $(W, \alpha)$, so we get $(W, \alpha \gamma)$ is
conjugate to $(W, \alpha \beta)$.  The conjugating element $x$, then, must
necessarily normalize $W$, and therefore, since $\alpha^G$ is
irreducible, $x \in W$.  Thus $\beta = \delta$.
\end{pf}

Now that we have proven Theorem 4.1, we are ready to prove the first main theorem of the introduction.  The following is a slight restatement of Theorem 1.1.  Essentially, we use the linear $\pi'$-special characters of the nucleus of the $\npi$ lift of $\varphi \in \ipig$ to construct lifts of $\varphi$.  We also point out that although this argument works equally well with subnormal nuclei and $\bpig$, the arguments in the following sections require normality instead of merely subnormality, so in order to relate these results to each other we need to use the normal nucleus and $\npig$.

\begin{theorem}  Let $G$ be a $\pi$-separable group, and let $\varphi \in \ipig$.  If $(W, \gamma)$ is the normal nucleus of the
unique lift of $\varphi$ in $\npig$, then for each linear character $\beta \in \mathcal{X}_{\pi'}(W)$, $(\gamma \beta)^G$ is a lift of $\varphi$, and if $\delta \in \mathcal{X}_{\pi'}(W)$ is linear, then $(\gamma \beta)^G = (\gamma \delta)^G$ if and only if $\beta = \delta$.  Thus $|W : W'|_{\pi'} \leq |L_{\varphi}|$.

\end{theorem}

\begin{pf}  Let $(W, \gamma)$ be the normal nucleus for
the unique character $\chi \in \npig$ such that $\chi^0 = \varphi$.  Let
$\lambda \in \irr(W/W')$ be $\pi'$-special.  Then by Theorem 4.1, $(\gamma \lambda)^G \in \irrg$.  Since $\lambda$ is a linear
$\pi'$-special character of $W$, then $(\gamma \lambda)^0 = \gamma^0$, and therefore $((\gamma \lambda)^G)^0 = (\gamma^0)^G = \varphi \in \ipig$.  Moreover, if $\lambda_1 \in \irr(W/W')$ is another linear $\pi'$-special character, then $(\gamma \lambda)^G = (\gamma \lambda_1)^G$ if and only if $\lambda = \lambda_1$.  Thus we have constructed an injective map $$\lambda \rightarrow (\gamma \lambda)^G$$ from the set of linear $\pi'$-special characters of $W$ into $L_{\varphi}$, and we therefore have $| W: W'|_{\pi'} \leq | L_{\varphi}| $.
\end{pf}

We now give an example to show that this lower bound may be strict.

\begin{example}  There exists a solvable group $G$ and a character $\varphi \in \ipig$ such that if $W$ is the normal nucleus of $\varphi$, then $1 = |W: W'|_{\pi'} < |L_{\varphi}|$.
\end{example}

Let $\Gamma$ be the nonabelian group of order 21.
Let $\Gamma$ act on $E$, an elementary abelian group of order $5^{21}$
such that every subgroup of $\Gamma$ is a stabilizer of some character
of $E$ (see \cite{counterexample} for details of this construction).
Let $G$ be the semidirect product of $\Gamma$ acting on $E$, and let $K \subseteq G$ be such that $E \subseteq K$ and $|G:K| = 3$.  Note $K \triangleleft G$.  Let $\pi = \{3, 5 \}$.  Choose a character $\alpha \in \irr(E)$ such that $G_{\alpha} = K$.  Now $\alpha$
is invariant in $K$, thus $\alpha$ must extend to a $\pi$-special character 
$\hat{\alpha} \in \irr(K)$, and $\beta = (\hat{\alpha})^G$ is
necessarily $\pi$-special.  Note that if $\beta^0 = \varphi \in \ipig$, then since $\beta$ is $\pi$-special, $W = G$, and
thus $| W: W '| _{\pi'} = 1$.  However, for each linear
$\pi'$-special character $\lambda \in \irr(K/E)$, we have that $\hat{\alpha} \lambda \in \irr(K)$ and $((\hat{\alpha} \lambda)^G)^0 = \varphi$.  If
$(\hat{\alpha}\lambda_1)^G =(\hat{\alpha}\lambda_2)^G$, note that
$\lambda_1$ and $\lambda_2$ must be conjugate in $G$, but all of the
elements of $G$ not in $K$ move $\hat{\alpha}$, thus this forces
$\lambda_1 = \lambda_2$.  Therefore $|L_{\varphi}|  \geq 7$, and in fact
it is easily seen that $| L_{\varphi}| = 7$.

\section{Vertices and correspondences}

In this section we discuss some correspondences between certain sets of characters that will be necessary to prove the upper bound stated in the introduction.  In particular, we will discuss Navarro's star map, which relates certain irreducible characters of an odd group $G$ with certain characters of a subgroup of $G$.  In the next section we will extend Navarro's result to provide a connection between lifts of a Brauer character of an odd group $G$ to certain Brauer characters of subgroups of $G$.

For a set of primes $\pi$, let $\Q_{\pi}$ denote the field obtained by adjoining all complex $n$th roots of unity of $\Q$ for all $\pi$-numbers $n$.  Recall (see \cite{bpi}, Corollary 12.1) that if $G$ is a $\pi$-separable group and if $2 \in \pi$ or $|G|$ is odd, and if $\chi \in \bpig$, then $\chi(g) \in \Q_{\pi}$ for every element $g \in G$.  The following lemma is well known, though the proof does not seem to appear in the literature.

\begin{lemma}  Let $G$ be a group of odd order, and let $\chi \in \irrg$.  If $\chi$ is $\pi$-factorable and $\chi^0 \in \ipig$, then $\chi$ has $\pi$-degree.
\end{lemma}

\begin{pf}  Suppose that $\chi$ factors as $\alpha \beta$, where $\alpha \in \mathcal{X}_{\pi}(G)$ and $\beta \in \mathcal{X}_{\pi'}(G)$.  Since $\chi^0 \in \ipig$, then $\beta^0 \in \ipig$ and thus the values of $\beta^0$ must be in $\Q_{\pi}$.  However, since $\beta$ is $\pi'$-special, then the values of $\beta$ must be in $\Q_{\pi'}$, and thus the values of $\beta^0$ must be in $\Q$, and therefore $\beta^0 = \bar{\beta}^0$.  Let $\eta \in \bpig$ be such that $\eta^0 = \beta^0$.  Thus $\bar{\eta}$ is the unique lift of $\bar{\beta}^0$, and therefore $\eta = \bar{\eta}$.  Since $G$ has odd order, this implies $\eta = 1_G$, and therefore $\beta$ is linear and $\chi = \alpha \beta$ has $\pi$-degree.
\end{pf}

One can show that $GL_2(3)$ is a counterexample to the above lemma if $|G|$ is not assumed to be odd.
 
\begin{corollary}  Suppose $G$ is a group of odd order, and suppose $\chi \in \irrg$ is a lift of $\varphi \in \ipig$.  Let $(W, \alpha \beta)$ be a normal nucleus for $\chi$.  Then $\beta$ is linear.
\end{corollary}

\begin{pf}  Note that since $((\alpha \beta)^G)^0 \in \ipig$, then $(\alpha \beta)^0 \in \ipi(W)$, and thus $\beta$ is linear by the above lemma.
\end{pf}

Recall (see \cite{oddcorrespondence}, for example) that if $p$ is a prime and $G$ is any finite group (we do not need to assume that $G$ is even $p$-solvable), and if $N \triangleleft G$ and $\theta \in \irr(N)$, then we say that $\chi \in \irr(G | \theta)$ is a relative defect zero character if $(\chi(1) / \theta(1))_p = |G : N|_p$.  We denote this by $\chi \in {\textup{rdz}}_p(G | \theta)$, or if $p$ is clear from the context, just ${\textup{rdz}}(G | \theta)$.  For instance, it is clear that if $N \subseteq H \subseteq G$, and if $\eta \in \irr(H | \theta)$ is such that $\eta^G = \chi$, then $\eta \in {\textup{rdz}}(H | \theta)$ if and only if $\chi \in {\textup{rdz}}(G | \theta)$.

The following lemma is in fact true without any assumptions on $G$ (see \cite{blocks}).  However, the proof given below seems easier than the proof in the general case and makes use of much of the theory we have already developed.

\begin{theorem}  Suppose that $G$ is a group of odd order and $Q$ is a normal $p$-subgroup of $G$.  Suppose that $\delta \in \irr(Q)$ is linear and invariant in $G$.  If $\chi \in {\textup{rdz}}(G | \theta)$, then $\chi^0 \in \ibrg$.
\end{theorem}

\begin{pf}  We prove this theorem by induction on the order of $G$.  Note that since $Q$ is a normal $p$-subgroup, then $(Q, \delta)$ is trivially a factored normal pair in $G$.  Thus we can choose a maximal factored normal pair $(N, \theta)$ such that $Q \subseteq N$ and $\theta$ lies between $\delta$ and $\chi$.

Suppose first that $\chi$ is not $p$-factorable.  Therefore $N < G$ and $G_{\theta} < G$.  Let $\eta \in \irr(G_{\theta} | \theta))$ be such that $\eta^G = \chi$.  Since $\chi \in {\textup {rdz}}(G | \delta)$, then $\eta \in {\textup{rdz}}(G_{\theta}| \delta)$, and the inductive hypothesis yields $\eta^0 \in \ibr(G_{\theta})$.  Let $(W, \xi)$ be a normal nucleus pair for for $\eta$ such that $N \subseteq W$ and $\theta$ lies under $\xi$.  Therefore $(W, \xi)$ is also a normal nucleus pair for $\chi$.  Note that $\xi^0$ must be in $\ibr(W)$.  By the previous corollary, $\xi$ has $p'$-degree, and thus $\xi_p$ must be linear.  Also note that since $\chi \in {\textup{rdz}}(G | \delta)$, then $\xi \in {\textup{rdz}}(W | \delta)$.  Therefore it must be the case that both $W/Q$ and $N/Q$ are $p'$-groups.

Note that $\xi_p$ must restrict to the linear $p$-special character $\theta_p \in \irr(N)$, which in turn must restrict to $\delta \in \irr(Q)$.  By the uniqueness of the $p$-special extension, we see that $G_{\theta_p} = G_{\delta} = G$, and thus $G_{\theta} = G_{\theta_{p'}} = G_{\theta^0}$.  By the Clifford correspondence, we therefore have that $(\eta^0)^G  \in \ibrg$, and since $(\eta^0)^G = \chi^0$, then we are done.

Otherwise, we may assume that $\chi$ is $p$-factorable.  Write $\chi = \alpha \beta$, where $\alpha \in {\mathcal{X}}_p(G)$ and $\beta \in {\mathcal{X}}_{p'}(G)$.  It is clearly enough to show that $\beta$ is linear.  Note that since $\chi \in {\textup{rdz}}(G | \delta)$ and $\delta$ is linear, then $\beta(1) = |G:Q|_p$.  Let $P \in {\textup{Syl}}_p(G)$, and set $\gamma = \beta_P \in \irr(P)$.  Since $Q$ is contained in the kernel of $\alpha$, then $\chi_Q = (\alpha \beta)_Q = \alpha_Q \beta_Q = e\beta_Q$ where $e$ is some positive integer.  Thus $\delta \in \irr(Q)$ lies under $\beta_P \in \irr(P)$.  Now $|P:Q| = \chi(1)_p = \beta(1)$, so $\delta^P = \beta_P$.  If $Q < P$, this contradicts the assumption that $\delta$ is invariant in $G$.  Therefore $Q = P$ and $\beta_P = \delta$, which is linear, and we are done.

\end{pf}

The following important result, from \cite{oddcorrespondence}, constructs a bijection between $\pi$-special characters and characters of certain subgroups.  We define, for a normal subgroup $N$ of $G$ and a Hall $\pi'$-subgroup $H$ of $G$, the set ${\mathcal{X}}_{\pi, H}(N)$ to be the set of $H$-invariant $\pi$-special characters of $N$.

\begin{theorem}  Let $G$ be a group of odd order and let $Q$ be a $p$-subgroup of $G$.  For every normal subgroup $N$ of $G$ such that $Q \cap N \in {\textup{Syl}}_p(N)$, there is a natural bijection $$^{\sim}: {\mathcal{X}}_{p', Q}(N) \rightarrow \irr(\nnq).$$  Moreover, suppose that $M \triangleleft G$ is such that $M \subseteq N$ and $\theta \in {\mathcal{X}}_{p', Q}(N)$ and $\eta \in {\mathcal{X}}_{p', Q}(M)$.  Then $\tilde{\theta}$ lies over $\tilde{\eta}$ if and only if $\theta$ lies over $\eta$.
\end{theorem}

Let ${\textup{Y}}_{p}(G) \sbs \ibrg$ denote the characters in $\ibrg$ with
$p'$-degree, and assume that the order of $G$ is odd.
Suppose $\chi$ is a $p'$-special character of $G$ such that $\chi^0 =          \varphi \in \ibrg$ and let $Q$ be a Sylow $p$-subgroup of $G$.  Since
${\textup{\textbf{N}}}_G(Q)/ Q$ is necessarily a $p'$-group, then
$(\tilde{\chi})^0$ is trivially in $\ibr({\textup{\textbf{N}}}_G(Q)/ Q)$.  In this case we abuse notation and say $\varphi \mapsto \tilde{\varphi}$ is a map from ${\textup{Y}}_{p}(G)$ to $\ibr({\textup{\textbf{N}}}_G(Q)/ Q)$.  The following result from \cite{oddweights} essentially extends this correspondence to a map from the Brauer characters of an odd group $G$ to irreducible Brauer characters of certain subgroups of $G$.  We use the notation $\varphi \in \ibr(G | Q)$ to signify that the vertex (see section 3) subgroup of $\varphi$ is $Q$.  

\begin{theorem}
Let $G$ be a group of odd order and let $\varphi \in \ibrg$.  Suppose $W$ is a subgroup of $G$ such that there exists a Brauer character $\alpha \in \ibr(W)$ of $p'$-degree such that $\alpha^G = \varphi$, and suppose $Q$ is a Sylow $p$-subgroup of $W$.  Then $\wt{\alpha} \in \ibr(\nwq)$ induces irreducibly to $(\wt{\alpha})^{\ngq} \in \ibr(\ngq)$.  Moreover, the map from $\ibr(G | Q)$ to $\ibr(\ngq | Q)$ given by $\varphi \rightarrow (\wt{\alpha})^{\ngq}$ is a well defined natural bijection.
\end{theorem}

Note that in theory, for a given Brauer character $\varphi$ of $G$, there are many subgroups $W$ which might have a Brauer character $\alpha$ that satisfy the hypotheses of the above theorem.  However, by Theorem 3.1, the subgroup $Q$ is unique up to conjugacy, and therefore part of the content of the above theorem is that the particular choice of $W$ and $\alpha$ does not affect the image of $\varphi$ in $\ibr(\ngq)$ under this map.

Now suppose that $\chi \in \irrg$ is a lift of the irreducible Brauer character $\varphi \in \ibrg$, and suppose that $G$ has odd order.  Note that by Corollary 5.2, if $(W, \gamma \delta)$ is a normal nucleus for $\chi$, then $\delta(1) = 1$ and thus $(\gamma \delta)^0 = \gamma^0$ is an irreducible Brauer character of $p'$-degree, and thus the pair $(W, \gamma^0)$ satisfies the hypotheses of the above theorem.  Moreover, by Theorem 3.1, any Sylow $p$-subgroup $Q$ of $W$ is conjugate to a vertex of $\varphi$.  The following theorem essentially extends the correspondence in the above theorem to lifts of $\varphi$, and shows that the images of distinct lifts of $\varphi$ must map to the same Brauer character of $\ngq$.

\begin{theorem}  Let $G$ be a group of odd order and suppose $\varphi \in \ibrg$.  Suppose $\chi, \psi \in \irrg$ are lifts of $\varphi$, and suppose that $\chi$ has normal nucleus $(W, \gamma)$ and $\psi$ has normal nucleus $(V, \epsilon)$, and suppose that a vertex subgroup $Q$ of $\varphi$ is contained in both $W$ and $V$.  Then $(\wt{\gamma^0})^{\ngq} \in \ibr(\ngq)$ and $(\wt{\gamma^0})^{\ngq} = (\wt{\epsilon^0})^{\ngq}$.
\end{theorem}

\begin{pf}  The first statement follows from Theorem 5.5, since $\gamma \in \irr(W)$ has $\pi$-degree and $\gamma^0 \in \ibr(W)$.  Since $\chi^0 = \psi^0 = \varphi$, then the injectivity in Theorem 5.5 implies that $(\wt{\gamma^0})^{\ngq} = (\wt{\epsilon^0})^{\ngq}$.
\end{pf}

In order to prove our upper bound on the number of lifts of the Brauer character $\varphi$, we will need to combine the above result with a result about Navarro's star map, which we now discuss.

Recall (see section 3) that if $G$ is $\pi$-separable with $\chi \in \irrg$ and if $\chi$ has normal nucleus $(W, \alpha \beta)$ and if $Q$ is any Hall $\pi'$-subgroup of $W$, then we say the pair $(Q, \beta_Q)$ is the vertex of $\chi$, and we denote this by ${\textup{vtx}}(\chi) = (Q, \beta_Q)$.  Note that $\beta_Q$ is necessarily irreducible.  For a $\pi'$-subgroup $Q$ and a character $\delta \in \irr(Q)$, recall that we use $\irr(G | Q, \delta)$ to denote all of the irreducible characters of $G$ that have $(Q, \delta)$ for the vertex.  Finally, we adopt the notation that $G_{\delta}$ is the stabilizer in $\ngq$ of the character $\delta \in \irr(Q)$.

The following theorem is, for our purposes, the key result from \cite{oddcorrespondence}.  In that paper, Navarro constructs, for a group of odd order $G$, a map $\chi \rightarrow \chi^*$ from $\irr(G | Q, \delta)$ to ${\textup{rdz}}(G_{\delta} | \delta)$.

\begin{theorem}  Suppose that $G$ is a group of odd order, and let $Q$ be a $p$-subgroup of $G$ and $\delta \in \irr(Q)$.  Then the map $\chi \rightarrow \chi^*$ is a natural injection from $\irr(G | Q, \delta)$ to ${\textup{rdz}}(G_{\delta} | \delta)$.
\end{theorem}

Before we discuss some of the specifics of the construction of the star map, we first prove a very useful corollary.

\begin{corollary}
Suppose that $G$ is a group of odd order, and suppose $\chi \in \irrg$ is such that $\chi^0 \in \ibrg$.  Moreover, suppose $\chi \in \irr(G | Q, \delta)$.  Then $(\chi^*)^0 \in \ibr(G_{\delta})$.
\end{corollary}

\begin{pf}
Since $\chi$ is a lift of some Brauer character of $G$, then by Corollary 5.2, the vertex character $\delta$ of $\chi$ is linear.  Thus Theorem 5.3 applied to $G_{\delta}$ implies that $(\chi^*)^0 \in \ibr(G_{\delta})$.
\end{pf}

For the full details of the construction of the star map, see \cite{oddcorrespondence}.  We will note, however, some important features of the construction that will be necessary in the next section.

\noindent \textbf{(1)}  If $\chi \in \irr(G| Q, \delta)$, then
$\chi^* \in \irr(G_{\delta})$ lies over $\delta$.

\medskip

\noindent \textbf{(2)}  Suppose $(W, \gamma)$ is a normal nucleus for the character $\chi \in \irr(G | Q, \delta)$ such that $(Q, \delta) \leq (W, \gamma)$, and suppose $(N, \theta)$ is a maximal factorable normal pair such that $(N, \theta) \leq (W, \gamma)$.  Note that since $Q$ is a Sylow $p$-subgroup of $W$, and $N \subseteq W$, then $Q \cap N$ is a Sylow $p$-subgroup of $N$.
Applying Theorem 5.4 to $N$ yields a $p'$-special character
$\widetilde{\theta_{p'}} \in \irr(\nnq)$.  Also, since $\theta_p \in {\mathcal{X}}_p(N)$ restricts irreducibly to $Q \cap N$, then $\theta_p$ restricts irreducibly to $(\theta_p)_{\nnq} \in {\mathcal{X}}_p(\nnq)$.  Let $\eta = \widetilde{\theta_{p'}} (\theta_p)_{\nnq} \in \irr(\nnq)$.  By the
construction of the star map, we have that $\eta = \widetilde{\theta_{p'}} (\theta_p)_{\nnq}$ lies under $\chi^*$.

\medskip

\noindent \textbf{(3)}  While keeping the same notation as in the
previous paragraph, we also let $T = G_{\theta}$.  It is shown in the
construction of the star map that $T_{\delta}$ is the stabilizer in
$G_{\delta}$ of $\eta \in \irr(\nnq)$.  Since $Q$ is a Sylow
$p$-subgroup of $W$, then we can again apply Theorem 5.4 to obtain a
$p'$-special character $\widetilde{\gamma_{p'}} \in \irr(\nwq)$.  Also,
since $|W: Q|$ is a $p'$-number, and $\gamma_p$ is $p$-special, then
$\gamma_p$ restricts irreducibly to $(\gamma_p)_{\nwq}$, and thus
$\widetilde{\gamma_{p'}} (\gamma_p)_{\nwq} \in \irr(\nwq)$.  It is shown
in the construction of the star map that $(\widetilde{\gamma_{p'}}             (\gamma_p)_{\nwq})^{T_{\delta}}$ is irreducible, and that
$(\widetilde{\gamma_{p'}} (\gamma_p)_{\nwq})^{T_{\delta}}$ is in fact
the Clifford correspondent for $\chi^*$ lying over $\eta \in                   \irr(\nwq)$.  Note then that $\chi^* = (\widetilde{\gamma_{p'}}                (\gamma_p)_{\nwq})^{G_{\delta}}$.

\medskip

\noindent \textbf{(4)}  Again, keeping the same notation as above,
note that $\gamma^T \in \irr(T | \theta)$ is the Clifford
correspondent for $\chi$ lying over $\theta$.  By the construction of the star map, $(\gamma^T)^* = (\widetilde{\gamma_{p'}} (\gamma_p)_{\nnq})^{T_{\delta}}  \in {\textup{rdz}}(T_{\delta} | \delta)$.

\section{The upper bound}

In this section we prove Theorem 1.2 of the introduction, which gives an upper bound on the number of lifts that a Brauer character of a group of odd order may have, in terms of the vertex subgroup $Q$.  Recall that Theorem 5.6 shows that distinct lifts of the same Brauer character $\varphi$ of $G$ map to the same Brauer character of $\ngq$, where $Q$ is the vertex subgroup of $\varphi$.  In this section we will modify Navarro's proof of the injectivity of the star map to show that distinct lifts of $\varphi$ with the same vertex character $\delta$ map to distinct Brauer characters of $G_{\delta}$, and we will then use this to show the upper bound on the number of lifts.

As in \cite{oddcorrespondence}, the following purely group theoretical fact will be used.

\begin{theorem}  Let $Q$ be a $p$-subgroup of $G$ and let $N$ and $M$ be normal subgroups of $G$.  Suppose that $Q \in {\textup{Syl}}_p(QM)$ and $Q \in {\textup{Syl}}_p(QN)$.  Then $\nnmq = \nnq \nmq$.
\end{theorem}

We now extend Navarro's argument to prove injectivity of the star map in \cite{oddcorrespondence} to prove injectivity of a related map.  Recall that Corollary 5.8 shows that if $\chi \in \irr(G|Q, \delta)$ is such that $\chi^0 \in \ibrg$, then $(\chi^*)^0 \in \ibr(G_{\delta})$.  Note that in the following theorem we are not assuming $\chi^0 = \mu^0$.    

\begin{theorem}  Let $G$ be a group of odd order and let $\pi = p'$.  Assume the characters $\chi$ and $\mu$ are in $\irr(G|Q, \delta)$, and suppose $\chi^0$ and $\mu^0$ are in $\ibrg$.  If $(\chi^*)^0 = (\mu^*)^0$, then $\chi = \mu$.
\end{theorem}

\begin{pf}  We will prove this theorem by induction on
  $|G|$.  For the reader's convenience, we will break the proof into
  several steps.

\noindent \textbf{Step 1}  If $\chi$ is $p$-factorable, then $\chi = \mu$.

\medskip

Since $\chi$ is $p$-factorable, then the pair $(G, \chi)$ is the
normal nucleus for $\chi$, and since $(Q, \delta)$ is a normal vertex
for $\chi$, then $Q$ is a full Sylow $p$-subgroup of $G$.  Since the
$p$-special factor $\chi_p$ of $\chi$ restricts to $\delta \in \irr(Q)$, and $|G : Q|$ is a $p'$-number, then $\delta$ extends to a
unique $p$-special character $\hat{\delta} \in {\mathcal{X}}_p(G_{\delta})$.  By the first fact about the star map, we see that $\chi^*$ and $\mu^*$ lie over $\delta$.

Since $\delta \in \irr(Q)$ extends to $\hat{\delta} \in                        {\mathcal{X}}_p(G_{\delta})$ and $|G_{\delta} : Q|$ is a $p'$-number,
then every irreducible character of $G_{\delta}$ lying over $\delta$
is $p$-factorable.  Thus $\chi^*$ and $\mu^*$ are both
$p$-factorable.

Note that since $\chi^0$ and $\mu^0$ are in $\ibrg$, then
Corollary 5.8 implies that $(\chi^*)^0$ and $(\mu^*)^0$ are in
$\ibr(G_{\delta})$.  We are assuming $(\chi^*)^0 = (\mu^*)^0$, and
  thus the $p'$-special factors of $\chi^*$ and $\mu^*$ are equal.
  Since the $p$-special factors of $\chi^*$ and $\mu^*$ are both equal
  to $\hat{\delta}$, then $\chi^* = \mu^*$.  By the injectivity of the
  star map, we have $\chi = \mu$.

\medskip

\noindent \textbf{Step 2}  From now on we may assume that neither
$\chi$ nor $\mu$ is $p$-factorable.  Let $(W, \gamma)$ be a normal
nucleus for $\chi$ such that $(Q, \delta) \leq (W, \gamma)$ and let
$(U, \rho)$ be a normal nucleus for $\mu$ such that $(Q, \delta) \leq (U, \rho)$.  Let $(N, \theta)$ be a maximal factorable normal pair such that
$(N, \theta) \leq (W, \gamma) \leq (G, \chi)$ and let $(M, \tau)$ be a maximal factorable normal pair such that $(M, \tau) \leq (U, \rho) \leq (G, \mu)$.  Then $(\gamma_p)_N = \theta_p$ and $\gamma_p(1) = 1$.
Similarly, $(\rho_p)_M = \tau_p$ and $\rho_p(1) = 1$.

\medskip

Since $(W, \gamma)$ is a normal nucleus for $\chi \in \irrg$, then
$\gamma^G = \chi$ and $\gamma$ is $p$-factorable.  Since $\chi^0 \in           \ibrg$, then by Corollary 5.2, $\gamma_p(1) = 1$.  Since $\gamma_p$ lies over $\theta_p \in \irr(N)$, then $(\gamma_p)_N = \theta_p$.  Similarly $\rho_p(1) = 1$ and $(\rho_p)_M = \tau_p$.

\medskip

\noindent \textbf{Step 3}  $(\theta_p)_{\nnq} = \widehat{\delta_{Q \cap N}}$, the unique $p$-special extension of $\delta_{Q \cap N}$ to $\nnq$.  Similarly, $(\tau_p)_{\nmq} = \widehat{\delta_{Q \cap M}}$, the unique $p$-special extension of $\delta_{Q \cap M}$ to $\nmq$.

\medskip

Note that $Q$ and $N$ are both contained in $W$, and thus since $Q$ is
a Sylow $p$-subgroup of $W$, then $Q \cap N$ is a Sylow $p$-subgroup
of $N$.  Since $\gamma_p \in \irr(W)$ is a linear $p$-special
character, then $\theta_p = (\gamma_p)_N$ is a linear $p$-special
character.  Thus $\theta_p \in {\mathcal{X}}_{p}(N)$ restricts to $(\theta_p)_{\nnq} \in {\mathcal{X}}_p(\nnq)$.  Since $\delta \in \irr(Q)$ also lies under the linear $p$-special character $\gamma_p$, then $\theta_p$ lies over
$\delta_{Q \cap N}$.  Therefore $(\theta_p)_{\nnq} = \widehat{\delta_{Q \cap N}}$, the unique $p$-special extension of $\delta_{Q \cap N}$ to $\nnq$.

Similarly, $(\tau_p)_{\nmq} = \widehat{\delta_{Q \cap M}}$, the unique
$p$-special extension of $\delta_{Q \cap M}$ to $\nmq$.

\medskip

\noindent \textbf{Step 4}  $\nnmq = \nnq \nmq$ and $\nnmq \subseteq G_{\delta}$.

\medskip

Recall $\gamma_p(1) = 1$ by Step 2, and since $\gamma_p \in \irr(W)$
lies over $\delta \in \irr(Q)$, then $\gamma_p$ is the unique
$p$-special extension of $\delta$ in $\irr(W)$.  Since $N \subseteq            W$, then any element of $N$ that normalizes $Q$ must stabilize
$\delta$.  Therefore $\nnq \subseteq G_{\delta}$, and similarly $\nmq          \subseteq G_{\delta}$.  By applying Theorem 6.1, we see that $\nnq \nmq        = \nnmq$, and thus $\nnmq \subseteq G_{\delta}$.

\medskip

\noindent \textbf{Step 5}  $Q \cap NM$ is a Sylow $p$-subgroup of $NM$
and $Q \cap NM$ is a normal Sylow $p$-subgroup of $\nnmq$.

\medskip

Since $Q$ is a Sylow $p$-subgroup of both $W$ and $U$, and $N \subseteq W$ and $M \subseteq U$, then $Q$ is a Sylow $p$-subgroup of both $QN$ and $QM$.  Thus $Q$ is a Sylow $p$-subgroup of $QMN$.  Thus $Q \cap NM$ is a Sylow $p$-subgroup of $NM$.  Since $Q \cap NM \subseteq \nnmq$, then $Q \cap NM$ is a normal Sylow $p$-subgroup of $\nnmq$.

\medskip

\noindent \textbf{Step 6}  There is an irreducible character $\alpha           \in {\mathcal{X}}_{p', Q}(NM)$ such that $\alpha$ lies over
$\theta_{p'} \in {\mathcal{X}}_{p'}(N)$ and $\alpha$ lies over
$\tau_{p'} \in {\mathcal{X}}_{p'}(M)$.

\medskip

Note that $\theta_{p'} \in {\mathcal{X}}_{p', Q \cap N}(N)$.  (Recall
that ${\mathcal{X}}_{p', Q \cap N}(N)$ denotes the $p'$-special
characters in $\irr(N)$ that are fixed by the $p$-subgroup $Q \cap N$.)
Therefore Theorem 5.4 applied to $N$ yields $\widetilde{\theta_{p'}} \in       \irr(\nnq / Q \cap N)$.  Note that $Q \cap N$ is a Sylow $p$-subgroup
of $N$, and thus $|\nnq : Q \cap N|$ is a $p'$-number, and thus
$\widetilde{\theta_{p'}} \in \irr(\nnq / Q \cap N)$ is necessarily
$p'$-special.  Similarly, Theorem 5.4 applied to $M$ yields
$\widetilde{\tau_{p'}} \in {\mathcal{X}}_{p'}(\nmq / Q \cap M)$.

Define $\eta \in \irr(\nnq)$ by $$\eta =
\widetilde{\theta_{p'}}(\theta_p)_{\nnq}.$$ By Step 3, we know that
$(\theta_p)_{\nnq} \in {\mathcal{X}}_p(\nnq)$, and thus $\eta$ is in
fact $p$-factorable and irreducible.  By the second fact about the
star map, we see that $\eta \in \irr(\nnq)$ lies under $\chi^*$.
Similarly, $\widetilde{\tau_{p'}} (\tau_p)_{\nmq} \in \irr(\nmq)$
lies under $\mu^*$.

By Step 2, we know that $(\theta_p)_{\nnq} \in {\mathcal{X}}_p(\nnq)$
is linear, and thus $\eta^0 = (\wt{\theta_{p'}})^0$.  Since
$(\chi^*)^0 \in \ibr(G_{\delta})$ and $\eta \in \irr(\nnq)$ lies
under $\chi$, then $\eta^0 = (\widetilde{\theta_{p'}})^0 \in                   \ibr(\nnq)$ lies under $(\chi^*)^0$.  Similarly,
$(\widetilde{\tau_{p'}})^0 \in \ibr(\nmq)$ lies under $(\mu^*)^0$,
and since we are assuming $(\chi^*)^0 = (\mu^*)^0$, then
$(\widetilde{\tau_{p'}})^0$ lies under $(\chi^*)^0$.

By Step 4, we know that $\nnmq \subseteq G_{\delta}$, and thus there
is a Brauer character $\epsilon \in \ibr(\nnmq)$ such that
$(\chi^*)^0 \in \ibr(G_{\delta})$ lies over $\epsilon$ and $\epsilon$ lies over $(\widetilde{\theta_{p'}})^0 \in \ibr(\nnq)$.  Choose an irreducible
Brauer character $\xi \in \ibr(\nmq)$ such that $\xi$ lies under
$\epsilon \in \ibr(\nnmq)$, and note that this implies that $\xi$ lies
under $(\chi^*)^0$.

Since both $\xi \in \ibr(\nmq)$ and $(\widetilde{\tau_{p'}})^0 \in            
\ibr(\nmq)$ lie under $(\chi^*)^0 \in \ibr(G_{\delta})$, and $\nmq \triangleleft G_{\delta}$, then $\xi$ is $G_{\delta}$-conjugate to
$(\widetilde{\tau_{p'}})^0$, and thus we can assume that $\xi =                (\widetilde{\tau_{p'}})^0$.

Note that $(\widetilde{\theta_{p'}})^0 \in \ibr(\nnq)$ and
$(\widetilde{\tau_{p'}})^0 \in \ibr(\nmq)$ both lie under $\epsilon            \in \ibr(\nnmq)$.  Thus, since $\nnmq = \nnq \nmq$, we see that $\epsilon \in \ibr(\nnmq)$ has $p'$-degree.  Therefore, there exists a character $\omega \in  {\mathcal{X}}_{p'}(\nnmq)$ such that $\omega^0 = \epsilon$.  By Step 5
we know that $Q \cap NM$ is a normal Sylow $p$-subgroup of $\nnmq$,
and thus $\omega \in \irr(\nnmq / Q \cap NM)$.  Therefore, by Theorem
5.4, there is a character $\alpha \in {\mathcal{X}}_{p', Q}(NM)$ such
that $\wt{\alpha} = \omega$.

Recall that $\epsilon \in \ibr(\nnmq)$ lies over
$(\widetilde{\theta_{p'}})^0 \in \ibr(\nnq)$, and thus $\omega \in {\mathcal{\
X}}_{p'}(\nnmq)$ lies over $\widetilde{\theta_{p'}} \in {\mathcal{X}}_{p'}(\nnq)$.  Again, by Theorem 5.4, we have that
$\alpha \in {\mathcal{X}}_{p', Q}(NM)$ lies over $\theta_{p'}$.
Similarly $\alpha$ lies over $\tau_{p'}$.

\medskip

\noindent \textbf{Step 7}  There exists a character $\hat{\delta} \in          {\mathcal{X}}_p(QMN)$ such that $\hat{\delta}$ is an extension of
$\delta \in \irr(Q)$ and $\hat{\delta}$ lies over both $\theta_p \in           {\mathcal{X}}_p(N)$ and $\tau_p \in {\mathcal{X}}_p(M)$.

\medskip

Since $\gamma_p \in {\mathcal{X}}_p(W)$ and $\rho_p \in                        {\mathcal{X}}_p(U)$ are linear and lie over $\delta \in \irr(Q)$, then
$\delta$ extends to $W$ and $U$.  Since $QN \subseteq W$ and $QM               \subseteq U$, then $\delta$ extends to $QN$ and $QM$.  Therefore
Theorem A of \cite{extension} applied to $QMN$ implies that $\delta$
extends to $QMN$.  Let $\hat{\delta} \in {\mathcal{X}}_p(QMN)$ be the unique
$p$-special extension of $\delta$ to $QMN$.  Since $(\gamma_p)_{QN}$
is a $p$-special extension of $\delta$, then $(\gamma_p)_{QN} =                \hat{\delta}_{QN}$.  Recall that $\gamma_p \in {\mathcal{X}}_p(W)$
lies over $\theta_p \in {\mathcal{X}}_p(N)$, and thus $\hat{\delta}$
lies over $\theta_p$.  Similarly $\hat{\delta}$ lies over $\tau_p$.

\medskip

\noindent \textbf{Step 8}  $N = M$ and $\theta = \tau$.

\medskip

Since $N \subseteq NM \triangleleft QNM$, there is a character
$\beta_1 \in {\mathcal{X}}_p(NM)$ such that $\beta_1$ lies under
$\hat{\delta} \in {\mathcal{X}}_p(QNM)$ and $\beta_1$ lies over
$\theta_p \in {\mathcal{X}}_p(N)$.  By Step 6, $\alpha \in {\mathcal{X}}_{p'}(NM)$ lies over $\theta_{p'} \in {\mathcal{X}}_{p'}(N)$, and thus $\alpha \beta_1$ is a $p$-factorable irreducible character of $NM$ lying over $\theta = \theta_p \theta_{p'} \in \irr(N)$.  Since $(N, \theta)$ is a maximal factorable pair,
then this forces $N = NM$ and $\theta = \alpha \beta_1$.  Note then
that $\beta_1 = (\hat{\delta})_N$, and thus $\theta = \alpha(\widehat{\delta})_N$.

Similarly, $M = MN$ and $\tau = \alpha (\hat{\delta})_M$.  Therefore
$M = N$ and $\theta = \alpha (\hat{\delta})_N = \tau$.

\medskip

\noindent \textbf{Step 9}  $\chi = \mu$.

\medskip

By Step 8, we now have that $N = M$ and $\theta = \tau$.  Set $T =             G_{\theta} = G_{\tau}$.  Notice that since $(N, \theta)$ is a maximal factorable normal pair and $N < G$, then we must have that $T < G$.

As in Step 6, we define $\eta \in \irr(\nnq)$ by $\eta =                       \widetilde{\theta_{p'}}(\theta_p)_{\nnq}$.  Then by the third fact
about the star map, we have that $T_{\delta} = (G_{\delta})_{\eta}$
and $((\gamma_p)_{\nwq} \wt{\gamma_{p'}})^{T_{\delta}} \in                     \irr(T_{\delta})$ is the Clifford correspondent for $\chi^*$ lying
over $\eta \in \irr(\nnq)$.  Similarly, $((\rho_p)_{\nuq}                      \wt{\rho_{p'}})^{T_{\delta}} \in \irr(T_{\delta})$ is the Clifford correspondent for $\mu^*$ lying over $\eta \in \irr(\nnq)$.

Now $G_{\delta}$ normalizes $Q$, and thus $G_{\delta}$ normalizes $Q           \cap N$ and $\nnq$.  Thus $G_{\delta}$ fixes $\delta_{Q \cap N}$ and
therefore by Step 3, $G_{\delta}$ stabilizes $\eta_p =                         (\theta_p)_{\nnq}$, the unique $p$-special extension of $\delta_{Q             \cap N}$ to $\nnq$.

Recall that by the third fact about the star map, we know that
$T_{\delta} = (G_{\delta})_{\eta}$.  Since $\eta_p$ is invariant in
$G_{\delta}$, then $T_{\delta}$ is the stabilizer in $G_{\delta}$ of
$\eta_{p'} = \wt{\theta_{p'}} \in {\mathcal{X}}_{p'}(\nnq)$ and thus
$T_{\delta} = (G_{\delta})_{\eta_{p'}} = (G_{\delta})_{\eta^0}$.  Note
that since $\eta \in \irr(\nnq)$ lies under $\chi^* \in \irr(G_{\delta})$, then $\eta^0 = (\theta_{p'})^0 \in \ibr(\nnq)$ lies under $(\chi^*)^0 = (\mu^*)^0  \in \ibr(G_{\delta})$.  Since $(\chi^*)^0 \in \ibr(G_{\delta})$
and $((\gamma_p)_{\nwq} \wt{\gamma_{p'}})^{T_{\delta}} \in                     \irr(T_{\delta})$ induces to $\chi^*$ by the third fact about the star
map, then $(((\gamma_p)_{\nwq} \wt{\gamma_{p'}})^{T_{\delta}})^0 \in           \ibr(T_{\delta})$.  Note that since $((\gamma_p)_{\nwq} \wt{\gamma_{p'}})^{T_{\delta}}$ lies over $\eta$, then $(((\gamma_p)_{\nwq}        
\wt{\gamma_{p'}})^{T_{\delta}})^0$ lies over $\eta^0$ and under
$(\chi^*)^0$, and thus $(((\gamma_p)_{\nwq}                                    \wt{\gamma_{p'}})^{T_{\delta}})^0$ is the Clifford correspondent for
$(\chi^*)^0 \in \ibr(G_{\delta})$ lying over $\eta^0 \in \ibr(\nnq)$.

Similarly $(((\rho_p)_{\nuq} \wt{\rho_{p'}})^{T_{\delta}})^0 \in               \ibr(T_{\delta})$ is the Clifford correspondent for $(\mu^*)^0 \in             \ibr(G_{\delta})$ lying over $\eta^0 \in \ibr(\nnq)$.  Since $(\chi^*)^0 = (\mu^*)^0$, then
$$(((\gamma_p)_{\nwq} \wt{\gamma_{p'}})^{T_{\delta}})^0 =
(((\rho_p)_{\nuq} \wt{\rho_{p'}})^{T_{\delta}})^0$$ by the uniqueness
of the Clifford correspondence.

By the fourth fact about the star map, $(\gamma^T)^* =                         ((\gamma_p)_{\nwq} \wt{\gamma_{p'}})^{T_{\delta}}$ and $(\rho^T)^* =           ((\rho_p)_{\nuq} \wt{\rho_{p'}})^{T_{\delta}}$.  Thus the above
paragraph implies $((\gamma^T)^*)^0 = ((\rho^T)^*)^0$.  Note that since
$(\gamma^T)^G = \chi$ and $\chi^0 \in \ibr(G)$, then $(\gamma^T)^0             \in \ibr(T)$.  Similarly $(\rho^T)^0 \in \ibr(T)$.  By Step 1 we can
assume that $T < G$, and therefore the inductive hypothesis yields
that $\gamma^T = \rho^T$.  Therefore $$\chi = \gamma^G = (\gamma^T)^G
= (\rho^T)^G = \rho^G = \mu,$$ and we are done.

\end{pf}

Before proving Theorem 1.2, we make some definitions.  If $\delta$ is a
linear character of a $p$-subgroup $Q$ of $G$, then we define $[               \delta ]$ to be the orbit of $\delta$ under the action of
${\bf{N}}_G(Q)$.  We also let ${\mathcal{O}}$ denote a set of
representatives of these orbits.  Finally, for a fixed character
$\varphi \in \ibr(G)$, we define $L_{\varphi}(\delta) =                       
L_{\varphi}([ \delta]) = L_{\varphi} \cap \irr(G | Q, \delta)$.
Also, recall that by Theorem 3.1, if $G$ has odd order and if the
character $\chi \in \irrg$ is a lift of the character $\varphi \in             \ibrg$, then the vertex subgroup is the same (up to conjugacy) as the
vertex subgroup of $\varphi$.

\begin{pf}  (of Theorem 1.2) Notice that by Theorem 3.1, if $\chi \in      
  \irrg$ is a lift of $\varphi \in \ibrg$, then the vertex subgroup of $\chi$ is (up to conjugation) $Q$, the vertex subgroup of $\varphi$.  Note that $L_{\varphi} = \bigcup\limits_{[\delta] \in \mathcal{O}} L_{\varphi}([\delta])$, and this is a disjoint union.  Thus we see that $|L_{\varphi}| = \sum\limits_{[\delta] \in \mathcal{O}} |L_{\varphi}([\delta])| $. 

We claim that $|L_{\varphi}([\delta])|  \leq | \ngq: G_{\delta}|$.  By Corollary 5.8, for each lift $\chi \in L_{\varphi}(\delta)$, we have
  $(\chi^*)^0$ is an irreducible Brauer character of $G_{\delta}$.
  Let $\chi \in \irrg$ and $\psi \in \irrg$ be distinct members of
  $L_{\varphi}(\delta)$, and suppose $(W, \gamma)$ is a normal nucleus for
  $\chi$ such that $(Q, \delta) \leq (W, \gamma)$ and $(U, \rho)$ is a
  normal nucleus for $\psi$ such that $(Q, \delta) \leq (U, \rho)$.
  Then Theorem 5.6 implies that $(\wt{\gamma^0})^{\ngq} =                
  (\wt{\rho^0})^{\ngq}$.  However, the third fact about the star map
  implies $(\chi^*)^0 = (\wt{\gamma^0})^{G_{\delta}}$ and
  $(\psi^*)^0 = (\wt{\rho^0})^{G_{\delta}}$.  Therefore
  $$((\chi^*)^0)^{\ngq} = ((\wt{\gamma^0})^{G_{\delta}})^{\ngq} =
  ((\wt{\rho^0})^{G_{\delta}})^{\ngq} = ((\psi^*)^0)^{\ngq}.$$
  However, the previous theorem states that $(\chi^*)^0$ is not equal
  to $(\psi^*)^0$.  Thus the images of the members of the set
  $L_{\varphi}(\delta)$ under the map given by $\chi \mapsto                    (\chi^*)^0$ must be distinct characters in $\ibr(G_{\delta})$, each
  of which induces to the same irreducible Brauer character of $\ngq$.
  Therefore, the size of the image of $L_{\varphi}([\delta])$ under
  the star map must be bounded above by $|\ngq : G_{\delta}|$. Since
  the star map is an injection, this proves that
  $|L_{\varphi}([\delta])| \leq | \ngq : G_{\delta}|$.

Therefore the sum becomes $$| L_{\varphi}|  \leq \sum\limits_{[\delta]
  \in \mathcal{O}} | \ngq : G_{\delta}| = | Q: Q'|$$ and we are done.
\end{pf}

The above result leads to some interesting questions.  For instance, is the requirement that the order of G is odd really necessary?  Certain aspects of the above proof seem to heavily depend on the oddness of $|G|$, but other arguments could be made.  One could probably begin by studying the case where $|G|$ is not necessarily odd, but $2 \in \pi'$, which seems to be the case where the lifts are more ``under control'' (see \cite{counterexample}).

Also, given the definition of the vertex $(Q, \delta)$ as above, one could ask which of the properties of a ``classical'' vertex extend to a vertex pair $(Q, \delta)$.  This might in turn lead to more precise statements about the number of lifts of a Brauer character.



\end{document}